\documentclass[11pt,a4paper,reqno]{amsart}
\usepackage{amsfonts,amsthm, amscd, epsfig, amsmath, amssymb,enumerate}
\usepackage{amstext}
\usepackage{xspace}
\usepackage{graphicx}
\usepackage{url}
\usepackage{enumerate}
\usepackage[all,2cell,ps]{xy}
\usepackage{xcolor}
\usepackage{subcaption}
\graphicspath{{./figure/}}
\usepackage{cite}
\usepackage{afterpage} 
\usepackage{faktor}

\newtheorem{itlemma}{Lemma}[section]
\newtheorem{itproposition}[itlemma]{Proposition}
\newtheorem{itfact}[itlemma]{Fact}
\newtheorem{theorem}[itlemma]{Theorem}
\newtheorem{itcorollary}[itlemma]{Corollary}
\newtheorem{itremark}[itlemma]{Remark}
\newtheorem{itremarks}[itlemma]{Remarks}
\newtheorem{itdefinition}[itlemma]{Definition}
\newtheorem{itexample}[itlemma]{Example}

\newenvironment{fact}{\begin{itfact}\rm}{\end{itfact}}
\newenvironment{claim}{\begin{itclaim}\rm}{\end{itclaim}}
\newenvironment{lemma}{\begin{itlemma}}{\end{itlemma}}
\newenvironment{remark}{\begin{itremark}\rm}{\end{itremark}}
\newenvironment{remarks}{\begin{itremarks} \rm}{\end{itremarks}}
\newenvironment{corollary}{\begin{itcorollary}}{\end{itcorollary}}
\newenvironment{proposition}{\begin{itproposition}}{\end{itproposition}}
\newenvironment{definition}{\begin{itdefinition}\rm}{\end{itdefinition}}
\newenvironment{example}{\begin{itexample}\rm}{\end{itexample}}
\newcommand{\be}[1]{\begin{equation}\label{#1}}
\newcommand{\ee}{\end{equation}}
\newcommand{\bl}[1]{\begin{lemma}\label{#1}}
\newcommand{\br}[1]{\begin{remark}\label{#1}}
\newcommand{\brs}[1]{\begin{remarks}\label{#1}}
\newcommand{\bt}[1]{\begin{theorem}\label{#1}}
\newcommand{\bd}[1]{\begin{definition}\label{#1}}
\newcommand{\bp}[1]{\begin{proposition}\label{#1}}
\newcommand{\bfact}[1]{\begin{fact}\label{#1}}
\newcommand{\bc}[1]{\begin{corollary}\label{#1}}
\newcommand{\bex}[1]{\begin{example}\label{#1}}
\newcommand{\ec}{\end{corollary}}
\newcommand{\efact}{\end{fact}}
\newcommand{\eex}{\end{example}}
\newcommand{\el}{\end{lemma}}
\newcommand{\er}{\end{remark}}
\newcommand{\ers}{\end{remarks}}
\newcommand{\et}{\end{theorem}}
\newcommand{\ed}{\end{definition}}
\newcommand{\ep}{\end{proposition}}
\newcommand{\epr}{\end{proof}}
\newcommand{\bpr}{\begin{proof}}
\newcommand{\bcl}[1]{\begin{claim}\label{#1}}
\newcommand{\ecl}{\end{claim}}

\newcommand{\ecs}{\end{corollary}}
\newcommand{\eers}{\end{exercise}}
\newcommand{\eexs}{\end{example}}
\newcommand{\eems}{\end{example}}
\newcommand{\els}{\end{lemma}}
\newcommand{\eles}{\end{lemmaex}}
\newcommand{\ets}{\end{theorem}}
\newcommand{\eds}{\end{definition}}
\newcommand{\eps}{\end{proposition}}

\newcommand{\bi}{\begin{itemize}}
\newcommand{\ei}{\end{itemize}}
\newcommand{\ben}{\begin{enumerate}}
\newcommand{\een}{\end{enumerate}}

\def\vbar{\mathchoice{\vrule height6.3ptdepth-.5ptwidth.8pt\kern-.8pt}
   {\vrule height6.3ptdepth-.5ptwidth.8pt\kern-.8pt}
   {\vrule height4.1ptdepth-.35ptwidth.6pt\kern-.6pt}
   {\vrule height3.1ptdepth-.25ptwidth.5pt\kern-.5pt}}
\def\fudge{\mathchoice{}{}{\mkern.5mu}{\mkern.8mu}}
\def\bbc#1#2{{\rm \mkern#2mu\vbar\mkern-#2mu#1}}
\def\bbb#1{{\rm I\mkern-3.5mu #1}}
\def\bba#1#2{{\rm #1\mkern-#2mu\fudge #1}}
\def\bb#1{{\count4=`#1 \advance\count4by-64 \ifcase\count4\or\bba A{11.5}\or
   \bbb B\or\bbc C{5}\or\bbb D\or\bbb E\or\bbb F \or\bbc G{5}\or\bbb H\or
   \bbb I\or\bbc J{3}\or\bbb K\or\bbb L \or\bbb M\or\bbb N\or\bbc O{5} \or
   \bbb P\or\bbc Q{5}\or\bbb R\or\bbc S{4.2}\or\bba T{10.5}\or\bbc U{5}\or
   \bba V{12}\or\bba W{16.5}\or\bba X{11}\or\bba Y{11.7}\or\bba Z{7.5}\fi}}

\def\1{{\bf 1}}

\def\ec{\`e }

\begin{document}

\title[]{Statistical test for an urn model with\\ random multidrawing
  and random addition}

\begin{abstract}
  We complete the study of the model introduced in
  \cite{cri-lou-min-multidrawing}.  It is a two-color urn model with
  multiple drawing and random (non-balanced) time-dependent
  reinforcement matrix. The number of sampled balls at each time-step
  is random. We identify the exact rates at which the number of balls
  of each color grows to $+\infty$ and define two strongly consistent
  estimators for the limiting reinforcement averages. Then we prove a
  Central Limit Theorem, which allows to design a statistical test for
  such averages.
\end{abstract}

\maketitle

    \centerline{Irene Crimaldi\footnote{IMT School for Advanced
        Studies Lucca, Piazza San Ponziano 6, 55100 Lucca, Italy,
        \url{irene.crimaldi@imtlucca.it}},
      Pierre-Yves~Louis\footnote{PAM UMR 02.102, Universit\'e
        Bourgogne Franche-Comt{\'e}, AgroSup Dijon, 1 esplanade Erasme,
        F-21000, Dijon, France,
        \url{pierre-yves.louis@agrosupdijon.fr}}\footnote{Institut de
        Math\'ematiques de Bourgogne, UMR 5584 CNRS, Universit\'e
        Bourgogne Franche-Comt{\'e}, F-21000, Dijon, France,
        \url{pierre-yves.louis@math.cnrs.fr}}, Ida
      G.~Minelli\footnote{Dipartimento di Ingegneria e Scienze
        dell'Informazione e Matematica, Universit\`a degli Studi
        dell'Aquila, Via Vetoio (Coppito 1), 67100 L'Aquila, Italy,
        \url{idagermana.minelli@univaq.it}}}

\bigskip
\noindent \textbf{Keywords.} Hypergeometric Urn; Multiple drawing urn;
P\'olya urn; Randomly reinforced urn; Central limit theorem; Stable
convergence; hypothesis testing; Population dynamics; Opinion
dynamics; Response-adaptive design.

\date{\today}

\section{Introduction}\label{intro}

\noindent An urn model can be described as follows: an urn contains
balls of different colors. At each time step a number $k$ of balls is
drawn and replaced into the urn together with other balls, whose color
depends on the composition of the sample.  It applies typically
whenever one wants to describe a \emph{(self-)reinforcement}
phenomenon in a stochastic evolution, which is the tendency to
increase the probability of an event in relation with the number of
times this event occurred in the past.\\ \indent We focus on urns with
balls of two colors, say $A$ and $B$.  The first urn model was from
P\'olya-Eggenberger \cite{egg-pol} and it consists in drawing
\emph{one} ball at each time and returning it back into the urn
together with $a$ balls of the same color. A variant of this model is
the Friedman urn \cite{fri}, where, at each time step, $a$ balls of
the color drawn and $b$ balls of the color not drawn are added into
the urn. A natural generalization consists in drawing more than one
ball and changing the replacement scheme at each time step. More
precisely, assume that the urn contains balls of color $A$ and $B$.
If $X_n$ (resp.,\ $k-X_n$) denotes the number of balls of color $A$
(resp.,\ color $B$) in the sample at time $n$, the number of balls of
color $A$ (resp.,\ $B$) added into the urn is $A_n X_n$
(resp.,\ $B_n(k-X_n)$), where the \emph{reinforcement factors} $A_n$
and $B_n$ may change over time. \\ \indent This class of models has
been studied by many authors in almost two decades and the usual
assumptions are that the number $k$ of sampled balls and the number
$(A_n+B_n)k$ of balls added into the urn at each step are constant
(e.g. \cite{chen-kuba,chen-wei-2005,Higueras2006,Idriss2018,
  Johnson2004,kuba2016classification,KubaMahmoud-balanced-affine-2017,mailler,mahmoud_2013_multisets}).\\ \indent
In more recent models, the replacement scheme is \emph{not balanced},
i.e., $A_n+B_n$ change in time and the reinforcement factors $A_n,$
and $B_n$ are \emph{random}. For example, in
\cite{AguechLasmarSelmi2019, AguechSelmi-unbalanced}, $(A_n)$ and
$(B_n)$ are two independent sequences of i.i.d. random
variables. Other authors have studied the case when the sample size,
denoted by $N_n$, is \emph{time-dependent, random} and \emph{possibly
  dependent of the past} and both colors are reinforced at each step
with the same i.i.d. reinforcement factors, i.e. $A_n=B_n$ for all $n$
(see \cite{perron, cri-ipergeom}, that include as a particular case
the model recently given in~\cite{Chen2020}).\\ \indent In
\cite{cri-lou-min-multidrawing}, the above models have been
generalized to the case when $A_n, B_n$ may be different and
eventually \emph{correlated}. The purpose of this article is to
complete the study of this last model by proving additional
theoretical results and providing a statistical tool for inference on
the reinforcement factors. \\ \indent Going into more detail, in
\cite{cri-lou-min-multidrawing}, the number of balls added into the
urn at each time-step is $(A_n+B_n)N_n$, of which $A_n X_n$ are
colored $A$ and $B_n(N_n-X_n)$ are colored $B$.  The pair of
reinforcement factors $[A_n, B_n]$ is assumed to be independent from
all the past until $n-1$ and from the composition of the sample at
time-step $n$. Thinking about possible applications of the model, this
assumption of independence means that the reinforcement deals with
some exogenous random factors, not related to those in the past and
external to the urn mechanism.  On the other hand, the joint
distribution of the pair $[A_n,B_n]$ may be whatever: the two random
factors may be correlated and their distribution may depend on
$n$. The main object of study is the asymptotic behavior of the urn
composition. For the proportion of color $A$ balls in the urn, almost
sure convergence results, as well as fluctuation results (through
Central Limit Theorems in the sense of stable convergence and of
almost sure conditional convergence), have been proven in
\cite{cri-lou-min-multidrawing}. Specifically, two settings have been
considered. If the factors $A_n$ and $B_n$ have the same mean (equal
reinforcement averages case), the limit proportion $Z\in [0,1]$ is
random without atoms and the proven CLTs provides asymptotic
confidence intervals for it.  In the case of unequal limit
reinforcement averages, i.e., when $E[A_n]\longrightarrow m_A,
E[B_n]\longrightarrow m_B$, with $m_A>m_B$, the proportion converges
almost surely to $1$: in the present paper we identify the exact rate
of convergence of such proportion towards its limit and we show that
the total number of color $A$ and color $B$ balls observed in the
samples up to time $n$ are quantities of order $n$ and $n^{m_B/m_A}$
respectively. Note that, in the applications, when convergence of the
proportion of color $A$ balls is slow (as can be seen in some of the
examples illustrated in \cite{cri-lou-min-multidrawing}), it could be
difficult to understand whether we are in the first or in the second
setting. Using the above mentioned result we are able to define two
strongly consistent estimators for $m_A$ and $m_B$ and we prove a CLT
for these estimators. We then use this result to design a statistical
test for the hypothesis $m_A=m_B$ versus $m_A>m_B$.\\ \indent Below we
discuss some possible applications of the model to {\em
  population/social dynamics} and {\em response-adaptive designs}.
\\ \indent Generally speaking, assume to have a population, where the
units can be of two types, say $A$ and $B$. The units already present
in the population are represented by the colored balls inside the
urn. At each time-step, a group (with random size $N_n$, possibly
dependent of the past) of units generates other units: each unit in
the group generates a certain number of new units of the same type. At
time-step $n$, the reinforcement factors $A_n$ and $B_n$ represent the
ability to generate new units for the two different sub-populations
(specific example are given below).  Indeed, it is reasonable to
assume the ability to give rise to new units to be specific for each
of the two sub-populations, depend on time and on exogenous random
factors. Moreover, the possible correlation between $A_n$ and $B_n$
models the competition or the cooperation between the two
sub-populations. If one of the reinforcement factors is eventually
larger in mean, then the associated sub-population will eventually
dominate. If both reinforcement factors are equal in mean then in the
long-run some random equilibrium takes place.  As specific examples,
one may think to bacterial populations or to the diffusion of genetic
variants of a virus, but also to examples in opinion dynamics: for
instance,
\begin{itemize}
\item (Elections) There are two candidates, $A$ and $B$, and people
  who have already (irreversibly) decided who they are going to vote
  for. They are color $A$ and color $B$ balls into the urn, that is
  the two sub-populations.  At each time-step, a random group $N_n$ of
  voters becomes ``active'' and try to persuade others to adopt their
  same choice. The ability to persuade may depend on several factors
  (intrinsic of the voters, good campaign or reputation of the
  candidate, etc.) and may change in time. Therefore, it is reasonable
  to assume that $A_n$ and $B_n$ are random and eventually
  correlated. Here we assume that the ``observable'' quantities are
  $N_n, X_n, A_n, B_n$, where $A_n$ and $B_n$ may be interpreted as
  trends associated with one or the other choice, which can be
  obtained through polls, study of the specific context, analysis of
  voting flows, etc.
 \item (Diffusion of a binary opinion through social networks) Imagine
   a connected population on the internet, where communities of agents
   share the same binary opinion on a given subject and grow in time
   through the addition of new ``followers'', who support one of these
   two opinions. In such communities some misperceptions may spread
   faster than others, regardless of their rationality.  This can
   generate what is called ``behavioral epidemics'', i.e., there is a
   group reinforcement in adopting a given opinion (or choice) $A$ due
   to a behavioral pattern perceived by the agents. In this case, the
   two different opinions ($A$ and its opposite $B$) identify the two
   sub-populations and we expect a strong correlation between $A_n$
   and $B_n$, where we interpret $B_n$ as the rate at which
   misperception $A$ is debunked with success. On the other hand,
   choices have a cost (psychological, material): this implies that
   the number $N_n$ of ``influencers'' who are imitated by others is
   random and may depend on the past, and the reinforcement factors
   $A_n, B_n$ may change along time. Again, in this case we assume
   that $N_n$, $A_n$ and $B_n$ are observable quantities: for example,
   on \emph{Twitter}, $A_n$ and $B_n$ may be desumed from the current
   number of ``retweets'' of comments or hashtags supporting the two
   opinions.
\end{itemize}
\indent This model can be also applied in the context of
response-adaptive designs (see \cite{hu-ros,ros,ros-lac} for
exhaustive reviews and see \cite{mayflournoy} for an urn model without
multidrawing, applied to clinical trials, which is included in the
model considered in the present paper as the very special case when
the sample-size $N_n$ is equal to $1$, $E[A_n]=m_A$ and $E[B_n]=m_B$
for each $n$). In such a design, units enter the experiment
sequentially and are allocated randomly to one out of two
``treatments'', $A$ or $B$, according to a rule that depends on the
previous allocations and the previous observed responses.  The
experimenter has two simultaneous goals: 1) collecting evidence to
determine the superior treatment and 2) biasing the allocations
towards the better treatment, in order to reduce the proportion of
units in the experiment that receive the inferior treatment. The
design driven by the model considered in this paper satisfies this
twice requirement when we take the reinforcement factors $A_n$ and
$B_n$ of the model as representative of the goodness of the two
treatments $A$ and $B$. Indeed, if one treatment is better than the
other in average, the design will allocate the units to the best
treatment with a proportion that converges almost surely to $1$. When
the two treatments are equivalent, the design allocates the proportion
of units with a random limit (see Remark \ref{rem-design}). \\ \indent
As an example of application in the industry sector, think to a firm
that has to select the size of its production for two different kinds
of products, say product $A$ and product $B$. At the same time, the
firm wants to get information on which of the two products is more in
demand by potential buyers. The firm may perform the following design
associated to our model.  At each time-step $n$, the total of the
production of the firm is $N_n$. The firm decides to produce $X_n$
products of type-$A$ and $N_n-X_n$ products of type-$B$, according to
a sampling without replacement of $N_n$ balls from an urn following
the dynamics here described. The factors $A_n$ and $B_n$, which will
affect the decision about the production at time-step $n+1$, are
related to a market survey on the two products, conducted at time-step
$n$. The market survey can be assumed to be independent of the
production of the firm and of the previous surveys, and this justifies
the hypothesis of independence for $[A_n, B_n]$ required in the
model. \\[10pt]

\indent The paper is organized as follows. In Section~\ref{model} we
recall the model and set the assumptions. In Section~\ref{main} we
state and prove the main theoretical results. Section~\ref{test}
provides the test on the reinforcement averages. Finally, in
Appendix~\ref{app-technical} we prove some technical results, and, in
Appendix~\ref{app-auxiliary} and in Appendix~\ref{app-stable-conv}, we
recall some auxiliary results and the definition of stable
convergence.


\section{Model specification}\label{model}

We consider the model introduced in
\cite{cri-lou-min-multidrawing}. An urn initially contains
$a\in{\mathbb N}\setminus\{0\}$ balls of color A and $b\in{\mathbb
  N}\setminus\{0\}$ balls of color B. At each discrete time $n\geq 1$,
we simultaneously (\textit{i.e.} without replacement) draw a random
number $N_n$ of balls. Let $X_n$ be the number of extracted balls of
color A.  Then we return the extracted balls into the urn together
with $A_nX_n$ balls of color A and $B_n(N_n-X_n)$ balls of color
B. More precisely, we take a probability space
$(\Omega,\mathcal{A},P)$ and some random variables
$N_n,\,X_n,\,A_n,\,B_n$ defined on it and such that, for each $n\geq
1$, we have:

\begin{itemize}
\item[(A1)] the conditional distribution of the random variable $N_n$ given 
$$[N_1,X_1, A_1,B_1\dots,N_{n-1}, X_{n-1}, A_{n-1},B_{n-1}]$$ is
  concentrated on $\{1,\dots,S_{n-1}\}$ where $S_{n-1}$ is the total
  number of balls in the urn at time~$n-1$, that is 
\begin{equation}
S_{n-1}=a+b+\sum_{j=1}^{n-1} A_j X_j+\sum_{j=1}^{n-1}B_j(N_j-X_j);
\end{equation}

\item[(A2)] the conditional distribution of the random variable $X_n$ given 
$$[N_1,X_1, A_1,,B_1\dots,N_{n-1}, X_{n-1}, A_{n-1},B_{n-1}, N_n]$$ is
  hypergeometric with parameters $N_n,\, S_{n-1}$ and $H_{n-1}$, 
  where $H_{n-1}$ is the total number of balls of color~A at time~$n-1$, that is 
\begin{equation}
H_{n-1}=a+\sum_{j=1}^{n-1} A_j X_j;
\end{equation}

\item[(A3)] the random vector $[A_{n},B_n]$ takes values in ${\mathbb
  N}\setminus\{0\}\times{\mathbb N}\setminus\{0\}$ and it is
  independent of
  $$[N_1,X_1,A_1,B_1,\ldots,N_{n-1}, X_{n-1},A_{n-1},B_{n-1}, N_n, X_n]\,.$$
  \end{itemize}
\indent According to the above notation, the random variable $X_n$
corresponds to the number of balls having the color A in a random
sample without replacement of size $N_n$ from an urn with $H_{n-1}$
balls of color A and $K_{n-1}=(S_{n-1}-H_{n-1})$ balls of color B.
The reinforcement rule is of the ``multiplicative'' type: indeed, at
$n\geq 1$, we add to the urn $A_nX_n$ balls of color $A$ and
$B_n(N_n-X_n)$ balls of color $B$. Therefore, the total number of
balls added to the urn, that is $A_nX_n+B_n(N_n-X_n)$, is random and
depends on $n$.\\
\indent Note that we do not specify the conditional distribution of
the random variable $N_n$ (the sample size) given the past
$$[N_1,X_1, A_1,B_1\dots,N_{n-1}, X_{n-1}, A_{n-1},B_{n-1}]$$ nor the
distribution of $[A_n,B_n]$ (the random reinforcement factors $A_n$
and $B_n$ may have different distributions, they may be correlated and
their joint and marginal distributions may vary with $n$). Several
examples can be found in \cite{cri-lou-min-multidrawing}.\\ \indent It
is worthwhile to remark that this model includes the Hypergeometric
Randomly Reinforced Urn studied in~\cite{perron,crimaldi2016} (take
$A_n=B_n$ for all $n$), which in turn includes the model recently
given in \cite{Chen2020}. In particular, two special cases are the
classical P\'olya urn (the case with $N_n=1$ and $A_n=B_n=k\in
\mathbb{N}\setminus\{0\}$ for each $n$) and the $2$-colors randomly
reinforced urn with the reinforcements for the two colors equal or
different in mean (the case with $N_n=1$ for each $n$ and $[A_n, B_n]$
arbitrarily random in ${\mathbb N}\setminus\{0\}\times{\mathbb
  N}\setminus\{0\}$).  Moreover, as told in Section \ref{intro},
previous literature (we refer to the quoted papers in
Sec.~\ref{intro}) deals with the case when the sample size $N_n$ is a
fixed constant, not depending on $n$, and/or the balanced case
(constant number of balls added to the urn each time).  \\

\indent We set $Z_n$ equal to the proportion of balls of color A in
the urn immediately after the $n$th update, that is $Z_0=a/(a+b)$ and
\begin{gather*}
  Z_n=\frac{H_n}{S_n}\quad\hbox{for } n\geq 1.
\end{gather*}
Moreover we set 
\begin{equation*}
\mathcal{F}_0=\{\emptyset, \Omega \},\quad\mathcal{F}_n=
\sigma\bigl(N_1,X_1,A_1,B_1,\ldots,N_n, X_n,A_n,B_n\bigr)
\quad\hbox{for } n\geq 1\,,
\end{equation*}
and 
\begin{equation*}
\mathcal{G}_n=\mathcal{F}_n\vee\sigma(N_{n+1}),\quad
\mathcal{H}_n=\mathcal{G}_n\vee\sigma(A_{n+1},B_{n+1})\quad\hbox{for } n\geq 0.
\end{equation*}

By the above assumptions and notation, we have
\begin{equation}\label{sp-cond-rinforzi}
E[A_{n+1}\,|\,\mathcal{G}_n]=E[A_{n+1}]\,,\qquad
E[B_{n+1}\,|\,\mathcal{G}_n]=E[B_{n+1}]
\end{equation}
and
\begin{equation}\label{sp-cond}
  \begin{split}
  E[X_{n+1}\,|\,\mathcal{H}_n]&=E[X_{n+1}\,|\,\mathcal{G}_n]=N_{n+1}Z_n\,,\\
  E[N_{n+1}-X_{n+1}\,|\,\mathcal{H}_n]&=E[N_{n+1}-X_{n+1}\,|\,\mathcal{G}_n]=
  N_{n+1}(1-Z_n)\,.
  \end{split}
\end{equation}
Finally, we set $\mathcal{X}_{n}=\{0\vee [N_{n}-(S_{n-1}-H_{n-1})],\dots,N_{n}\wedge
H_{n-1}\}$ and, for each $k\in \mathcal{X}_{n}$,   
\begin{equation}\label{hypergeometric-distr}
p_{n,k}=p_k(N_{n}, S_{n-1}, H_{n-1})=
\frac{\binom{H_{n-1}}{k}\binom{S_{n-1}-H_{n-1}}{N_{n}-k}}{\binom{S_{n-1}}{N_{n}}}\,.
\end{equation}

\indent For all $n$, set the reinforcement averages/means
$m_{A,n}=E[A_n]\geq 1$ and $m_{B,n}=E[B_n]\geq 1$.  We will assume
that the two sequences $(m_{A,n})_n$ and $(m_{B,n})_n$ respectively
converge to $m_A<+\infty$ and $m_B<+\infty$.  In the present paper, we
will consider the following cases:
\\

\indent {\sc Case $m_A>m_B$}: We have $m_A>m_B$ with 
\begin{equation}\label{cond-medie}
\sum_n \frac{|m_{A,n+1}m_B-m_A m_{B,n+1}|}{n}<+\infty;
\end{equation}

\indent {\sc Case $m_A=m_B$:} We have $m_{A,n}=m_{B,n}=m_n$ for each
  $n$ and so $m_A=m_B=m\in [1,+\infty)$.
  \\
  
\noindent Note that assumption \eqref{cond-medie} is trivially satisfied when
$|m_{A,n}-m_A|=O(n^{-\epsilon_A})$ and
$|m_{B,n}-m_B|=O(n^{-\epsilon_B})$, with $\epsilon_A>0$ and
$\epsilon_B>0$. In particular, it is verified when $m_{A,n}=m_{A}$ and
$m_{B,n}=m_B$ for each $n$ (constant average values).\\ \indent In the
sequel, we will refer to the above two cases as ``case $m_A>m_B$'' and 
``case $m_A=m_B$'', respectively. Therefore, {\em condition
  \eqref{cond-medie} will be tacitly assumed} when we are in case 
$m_A>m_B$. Moreover, for simplicity, throughout the paper we will
assume
$$
A_n\vee B_n\vee N_n\leq C\qquad\mbox{ for some (integer) constant } C.
$$ As mentioned in \cite{cri-lou-min-multidrawing}, sometimes this
assumption can be removed or replaced by an assumption of uniform
integrability, but we will not focus in the present work on this
possibility.


\section{Main results}\label{main}

In \cite{cri-lou-min-multidrawing}, without assumption
\eqref{cond-medie}, we proved that $H_n$ and $K_n=S_n-H_n$ go almost
surely to $+\infty$ and that, under the condition
\begin{equation}\label{cond-N}
E[N_{n+1}|\mathcal{F}_n]\stackrel{a.s.}\longrightarrow N,
\end{equation}
where $N$ is a (finite, strictly positive) random variable,
we have 
$$
\frac{H_n}{n}\stackrel{a.s.}\longrightarrow m_A N Z,
\qquad\mbox{and}\qquad
\frac{K_n}{n}\stackrel{a.s.}\longrightarrow m_B N(1-Z),
$$ where $Z=1$ in the case $m_A>m_B$ and a non-atomic random variable
with values in $[0,1]$ in the case $m_A=m_B$.  Therefore, $H_n$ and
$S_n$ grow like $n$ in both cases, while $K_n$ grows as $n$ only in
the case $m_A=m_B$. For the case $m_A>m_B$, we obtained $K_n =
o(n^\theta)$, that is $1-Z_n = o(n^{-(1-\theta)})$, for all
$\theta\in(m_B/m_A,1)$. With the following result, we prove here that,
under the additional assumption \eqref{cond-medie}, $K_n$ goes to
$+\infty$ exactly as $n^{m_B/m_A}$.

\begin{theorem}\label{rate-th} (Rate of $K_n$ in case $m_A>m_B$)\\
  Suppose to be in case $m_A>m_B$ and assume \eqref{cond-N}.  Then:
\begin{itemize}
\item $(K_n)_n$ increases to $+\infty$ as $n^{m_B/m_A}$;
\item $n^{1-m_B/m_A}(1-Z_n)\stackrel{a.s.}\longrightarrow
  \widetilde{Z}$, where $\widetilde{Z}$ is a random
  variable taking values in $(0,+\infty)$;
\item setting $N_{B,n}=\sum_{j=1}^n (N_j-X_j)$, we have
  $\frac{N_{B,n}}{n^{m_B/m_A}}\stackrel{a.s.}\longrightarrow
  \frac{m_A}{m_B}N\widetilde{Z}$.
  \end{itemize}
\end{theorem}

\begin{proof}
  The main part of the proof is the proof of Lemma \ref{lemma-limit}.
  Indeed, given this result, that is the fact that $H_n/K_n^{m_A/m_B}$
  converges almost surely to a random variable $Y$ taking values in
  $(0,+\infty)$, we have
  $$
  \frac{K_n}{n^{m_B/m_A}}=
  \frac{K_n}{H_n^{m_B/m_A}}\left(\frac{H_n}{n}\right)^{m_B/m_A}
  \stackrel{a.s.}
  \longrightarrow \widetilde{Y}=Y^{-m_B/m_A}(m_AN)^{m_B/m_A}\in (0,+\infty).
  $$
  As a consequence, we have
  $$
  n^{1-m_B/m_A}(1-Z_n)=n^{1-m_B/m_A}\frac{K_n}{S_n}=
  \frac{K_n}{n^{m_B/m_A}}\frac{n}{S_n}
  \stackrel{a.s.}\longrightarrow \widetilde{Z}=\frac{\widetilde{Y}}{m_AN}
  =\frac{Y^{-m_B/m_A}}{(m_A N)^{1-m_B/m_A}}
  \in (0,+\infty).
  $$
  Finally, we observe that
  $$
  n^{1-m_B/m_A}E[N_{n+1}-X_{n+1}|\mathcal{F}_{n}]=
  E[N_{n+1}|\mathcal{F}_n]n^{1-m_B/m_A}(1-Z_n)
  \stackrel{a.s.}\longrightarrow N\widetilde{Z}
  $$ and so, by Lemma \ref{lemma-app} (applied with
  $Y_j=j^{1-m_B/m_A}(N_j-X_j)$ and $\alpha_j=j^{1-m_B/m_A}$ and
  $\beta_n=n^{m_B/m_A}$), we get
  $$
  \frac{N_{B,n}}{n^{m_B/m_A}}=\frac{1}{n^{m_B/m_A}}\sum_{j=1}^n(N_j-X_j)
  \stackrel{a.s.}\longrightarrow \frac{m_A}{m_B}N\widetilde{Z}.
  $$
  \end{proof}

\begin{corollary}\label{consistent-est} (Strongly consistent estimators)\\
  Assume condition \eqref{cond-N} and
\begin{equation}\label{cond-Q}
  E[N_{n+1}^2|\mathcal{F}_n]\stackrel{a.s.}\longrightarrow Q,
\end{equation}
where $Q$ is a (strictly positive finite) random variable. 
  In both cases $m_A=m_B$ and $m_A>m_B$, 
\begin{equation}\label{def-est-NQ}
  \widehat{\mu}_n=\frac{\sum_{j=1}^n N_j}{n},\qquad
  \widehat{q}_{N,n}=\frac{\sum_{j=1}^n N^2_j}{n}\,,
\end{equation}
are strongly consistent estimators of the random variables $N$ and
$Q$.  Moreover, setting $N_{A,n}=\sum_{j=1}^nX_j$ and
$N_{B,n}=\sum_{j=1}^n(N_j-X_j)$, the random variables
 \begin{equation}\label{def-est-average}
  \widehat{m}_{A,n}=\frac{\sum_{j=1}^n A_jX_j}{N_{A,n}},
  \qquad
  \widehat{m}_{B,n}=\frac{\sum_{j=1}^n B_j(N_j-X_j)}{N_{B,n}}
\end{equation}
 are strongly consistent estimators of $m_A$ and $m_B$, respectively.
 Further, assuming
 \begin{equation}\label{cond-var}
   q_{A,n}=E[A_n^2]\longrightarrow q_A,\qquad q_{B,n}=E[B_n^2]\longrightarrow q_B,
 \end{equation}  
 \begin{equation}\label{def-est-square}
  \widehat{q}_{A,n}=\frac{\sum_{j=1}^n
    A_j^2X_j}{N_{A,n}}\qquad\mbox{and}\qquad 
  \widehat{q}_{B,n}=\frac{\sum_{j=1}^n
    B_j^2(N_j-X_j)}{N_{B,n}}
\end{equation}
are strongly consistent estimators of $q_A$ and $q_B$. Finally, assuming
\begin{equation}\label{cond-cov}
  q_{AB,n}=E[A_nB_n]\longrightarrow q_{AB}
\end{equation}
and $P(N=1)<1$, we have
\begin{equation}\label{def-est-average-prod}
  \widehat{q}_{AB,n}=\frac{\sum_{j=1}^n
    A_jB_jX_j(N_j-X_j)}{\sum_{j=1}^nX_j(N_j-X_j)}
  \stackrel{a.s.}\longrightarrow q_{AB}
  \qquad\mbox{on } \{N>1\}.
\end{equation}
\end{corollary}

\begin{proof}
  For $\widehat{\mu}_n$ and $\widehat{q}_{N,n}$, it is enough to apply
  Lemma \ref{lemma-app} using the two assumptions \eqref{cond-N} and
  \eqref{cond-Q}.  For the other random variables, we have to employ
  Lemma \ref{lemma-app}, using the suitable rates given in
  \cite{cri-lou-min-multidrawing} and in the above Theorem
  \ref{rate-th}.  As an example, we show the proof only of the
  statement for $\widehat{q}_{AB,n}$.  To this regards, we observe
  that
  $$
  E[X_{n+1}(N_{n+1}-X_{n+1})|\mathcal{H}_n]=
  Z_n(1-Z_n)\frac{S_n}{S_n-1}(N_{n+1}^2-N_{n+1})
  $$
  and so
  $$
  n^{1-m_B/m_A}E[X_{n+1}(N_{n+1}-X_{n+1})|\mathcal{F}_n]
  \stackrel{a.s.}\longrightarrow
  \begin{cases}
    Z(1-Z)(Q-N)\quad&\mbox{in the case  } m_A=m_B\\
    \widetilde{Z}(Q-N)\quad&\mbox{in the case  } m_A>m_B.
  \end{cases}
  $$
 Therefore, Lemma \ref{lemma-app} implies
  $$
  \frac{\sum_{j=1}^nX_j(N_j-X_j)}{n^{m_B/m_A}}\stackrel{a.s.}\longrightarrow
  \begin{cases}
    Z(1-Z)(Q-N)\quad&\mbox{in the case  } m_A=m_B\\
    \frac{m_A}{m_B}\widetilde{Z}(Q-N)\quad&\mbox{in the case  } m_A>m_B,
  \end{cases}
  $$
  Similarly, since
  $E[A_{n+1}B_{n+1}X_{n+1}(N_{n+1}-X_{n+1})|\mathcal{F}_n]=
  q_{AB,n+1}E[X_{n+1}(N_{n+1}-X_{n+1})|\mathcal{F}_n]$, 
  we get
  $$
  \frac{\sum_{j=1}^nA_jB_jX_j(N_j-X_j)}{n^{m_B/m_A}}
  \stackrel{a.s.}\longrightarrow
  \begin{cases}
    q_{AB}Z(1-Z)(Q-N)\quad&\mbox{in the case  } m_A=m_B\\
    q_{AB}\frac{m_A}{m_B}\widetilde{Z}(Q-N)\quad&\mbox{in the case  } m_A>m_B,
  \end{cases}
  $$
  and so $\widehat{q}_{AB,n}\stackrel{a.s.}\longrightarrow q_{AB}$
  on $\{Q-N>0\}=\{N>1\}$.
\end{proof}

\begin{remark}\label{rem-design} \rm The proportion of balls of color $A$ along time is
$$ \frac{N_{A,n}}{\sum_{j=1}^n N_j}=\frac{\sum_{j=1}^n
    X_j}{n}\frac{n}{\sum_{j=1}^nN_j}.
  $$ By the above Corollary~\ref{consistent-est}, $\sum_{j=1}^nN_j/n$
  converges almost surely to $N$, while, with a similar argument,
  $\sum_{j=1}^n X_j/n$ converges almost surely to $Z$, because
  $E[X_{n+1}|\mathcal{H}_n]=N_{n+1}Z_n$ and so, under \eqref{cond-N},
  we have $E[X_{n+1}|\mathcal{F}_n]\stackrel{a.s.}\to NZ$. Therefore,
  the proportion $N_{A,n}/\sum_{j=1}^n N_j$ converges to $1$ in the
  case $m_A>m_B$ and to a non-atomic random variable $Z\in [0,1]$ in
  the case $m_A=m_B$. This is the formalization of the property
  required for response-adaptive designs, mentioned in the
  Introduction.
\end{remark}

\begin{theorem}\label{asymptotic-normality}
  (Asymptotic normality of the estimators)\\ With the same notation as
  in Corollary \ref{consistent-est}, assume conditions \eqref{cond-N},
  \eqref{cond-Q}, \eqref{cond-var} and \eqref{cond-cov} and set
  $$
  \sigma^2_A=q_A-m_A^2,\qquad\sigma^2_B=q_B-m_B^2\quad\mbox{and}\quad
  c_{AB}=q_{AB}-m_Am_B.$$ In both cases $m_A>m_B$ and $m_A=m_B$, we
  have
  $$
  \left[
    \sqrt{N_{A,n}}\left(\widehat{m}_{A,n}-m_{A}\right),
    \sqrt{N_{B,n}}\left(\widehat{m}_{B,n}-m_{B}\right)
    \right]
  \stackrel{stably}\longrightarrow \mathcal{N}(0,\Sigma),
  $$
  where
  $$\Sigma=
  \begin{pmatrix}
    \Sigma_{AA} &\Sigma_{AB}\\
    \Sigma_{AB} &\Sigma_{BB}
  \end{pmatrix}
  $$
  with
$$
\Sigma_{AA}=\begin{cases}
  \sigma^2_A Q/N\;&\mbox{when } m_A>m_B\\
  \sigma^2_A\left[(1-Z)+ZQ/N\right]\;&\mbox{when } m_A=m_B
\end{cases},
\quad
\Sigma_{BB}=  \begin{cases}
    \sigma^2_B\;&\mbox{when }m_A>m_B\\
    \sigma^2_B \left[Z+(1-Z)Q/N\right]\;&\mbox{when } m_A=m_B
\end{cases}
$$
and
$$
\Sigma_{AB}=
\begin{cases}
  0\;&\mbox{in the case  } m_A>m_B\\
  c_{AB}(Q/N-1)\;&\mbox{in the case  } m_A=m_B.
\end{cases}
$$
\end{theorem}
Note that the asymptotic covariance is always zero in the case
$m_A>m_B$ and, in the case $m_A=m_B$, it is equal to zero when
$c_{AB}=0$ or only on the event $\{Q=N\}=\{N=1\}$.

\begin{proof} We are going to apply
  Theorem \ref{thm:triangular}. To this purpose, we define the
  $2$-dimensional random vector $\mathbf{T}_{n,k}$ as
  $[\mathbf{T}_{n,k}]_1=\frac{X_k(A_k-m_{A,k})}{\sqrt{n}}$ and
  $[\mathbf{T}_{n,k}]_2=\frac{(N_k-X_k)(B_k-m_{B,k})}{\sqrt{n^{m_B/m_A}}}$.
  Moreover, set $\mathcal{G}_{n,k}=\mathcal{G}_k$.  We observe that
  $$ E[X_k(A_k-m_{A,k})|\mathcal{G}_{n,k-1}]=
  E[X_k(A_k-m_{A,k})|\mathcal{G}_{k-1}]=
  E[X_k|\mathcal{G}_{k-1}]E[A_k-m_{A,k}]=0$$ and, similarly,
  $E[(N_k-X_k)(B_k-m_{B,k})|\mathcal{G}_{k-1}]=0$.  Therefore, for any
  fixed $n$, $(\mathbf{T}_{n,k})_{1\leq k\leq n}$ is a martingale
  difference array with respect to the filtration
  $(\mathcal{G}_{n,k})_{k\geq 0}$, which satisfies condition (c1) of
  Theorem \ref{thm:triangular}.  Moreover, we have
  \begin{equation*}
    \sum_{k=1}^n \mathbf{T}_{n,k}\mathbf{T}_{n,k}^{\top}=
    \begin{pmatrix}
      \frac{\sum_{k=1}^nX_k^2(A_k-m_{A,k})^2}{n}
      &\frac{\sum_{k=1}^NX_k(N_k-X_k)(A_k-m_{A,k})(B_k-m_{B,k})}
           {n^{m_B/m_A}n^{1/2(1-m_B/m_A)}}\\
           \frac{\sum_{k=1}^n X_k(N_k-X_k)(A_k-m_{A,k})(B_k-m_{B,k})}
                {n^{m_B/m_A}n^{1/2(1-m_B/m_A)}}
      &\frac{\sum_{k=1}^n (N_k-X_k)^2(B_k-m_{B,k})^2}{n^{m_B/m_A}}    
    \end{pmatrix}.
    \end{equation*}
  Now, we set $\sigma^2_{A,k}=q_{A,k}-m_{A,k}^2=E[A_k^2]-E[A_k]^2$ and we
  observe that
  \begin{equation*}
    \begin{split}
  &E[X_k^2(A_k-m_{A,k})^2|\mathcal{F}_{k-1}]=E[X_k^2|\mathcal{F}_{k-1}]\sigma^2_{A,k}
  =\\
  &\sigma^2_{A,k}\left\{
  E[N_k|\mathcal{F}_{k-1}]Z_{k-1}(1-Z_{k-1})
  \frac{S_{k-1}-E[N_k|\mathcal{F}_{k-1}]}{S_{k-1}-1}+
  E[N_k^2|\mathcal{F}_{k-1}]Z_{k-1}^2\right\}\\
&  \stackrel{a.s.}\longrightarrow \widetilde{\Sigma}_{AA}=
  \begin{cases}
\sigma^2_A Q\;&\mbox{when } m_A>m_B\\
  \sigma^2_A\left[Z(1-Z)N+Z^2Q\right]\;&\mbox{when } m_A=m_B
  \end{cases}
    \end{split}
  \end{equation*}
  and so, by Lemma \ref{lemma-app} (applied to
  $Y_j=X_j^2(A_j-m_{A,j})^2$, $\alpha_j=1$ and $\beta_n=n$), we have
  $$
  \sum_{k=1}^nX_k^2(A_k-m_{A,k})^2/n\stackrel{a.s.}\longrightarrow
  \widetilde{\Sigma}_{AA}.$$
  Similarly, setting  $\sigma^2_{B,k}=q_{B,k}-m_{B,k}^2=E[B_k^2]-E[B_k]^2$, we have
 \begin{equation*}
    \begin{split}
  &k^{1-m_B/m_A}E[(N_k-X_k)^2(B_k-m_{B,k})^2|\mathcal{F}_{k-1}]=\\
  &\sigma^2_{B,k}\left\{
  E[N_k|\mathcal{F}_{k-1}]k^{1-m_B/m_A}(1-Z_{k-1})Z_{k-1}
  \frac{S_{k-1}-E[N_k|\mathcal{F}_{k-1}]}{S_{k-1}-1}+
  E[N_k^2|\mathcal{F}_{k-1}]k^{1-m_B/m_A}(1-Z_{k-1})^2\right\}\\
 & \stackrel{a.s.}\longrightarrow
\widetilde{\widetilde{\Sigma}}_{BB}=  \begin{cases}
    \sigma^2_B N\widetilde{Z}\quad &\mbox{when } m_A>m_B\\
    \sigma^2_B \left[Z(1-Z)N+(1-Z)^2Q\right]\quad&\mbox{when } m_A=m_B
    \end{cases}
    \end{split}
  \end{equation*} 
  and so, by Lemma \ref{lemma-app} (applied to
  $Y_j=j^{1-m_B/m_A}(N_j-X_j)^2(B_j-m_{B,j})^2$,
  $\alpha_j=j^{1-m_B/m_A}$ and $\beta_n=n^{m_B/m_A}$), we have
  $\sum_{k=1}^n(N_k-X_k)^2(B_k-m_{B,k})^2/n^{m_B/m_A}
  \stackrel{a.s.}\longrightarrow
  \widetilde{\Sigma}_{BB}=\frac{m_A}{m_B}\widetilde{\widetilde{\Sigma}}_{BB}$.
  Finally, looking at the computations reported in the proof of
  Theorem \ref{consistent-est}, we have
  \begin{equation*}
    \begin{split}
  &k^{1-m_B/m_A}E[X_k(N_k-X_k)(A_k-m_{A,k})(B_k-m_{B,k})|\mathcal{F}_{k-1}]=
  c_{AB,k}k^{1-m_B/m_A}E[X_k(N_k-X_k)|\mathcal{F}_{k-1}]\\
  &\stackrel{a.s.}\longrightarrow
      \widetilde{\widetilde{\Sigma}}_{AB}=
 \begin{cases}
    c_{AB}Z(1-Z)(Q-N)\quad&\mbox{in the case  } m_A=m_B\\
    c_{AB}\widetilde{Z}(Q-N)\quad&\mbox{in the case  } m_A>m_B.
  \end{cases}
    \end{split}
  \end{equation*}
  Therefore, by Lemma \ref{lemma-app} (applied to
  $Y_j=j^{1-m_B/m_A}E[X_j(N_j-X_j)(A_j-m_{A,j})(B_j-m_{B,j})$,
    $\alpha_j=j^{1-m_B/m_A}$ and $\beta_n=n^{m_B/m_A}$), we have
    $$
\sum_{k=1}^n
  X_k(N_k-X_k)(A_k-m_{A,k})(B_k-m_{B,k})/n^{m_B/m_A}
  \stackrel{a.s.}\longrightarrow \frac{m_A}{m_B}
  \widetilde{\widetilde{\Sigma}}_{AB}
  $$
  and so, setting $\widetilde{\Sigma}_{AB}$ equal to
  $\frac{m_A}{m_B} \widetilde{\widetilde{\Sigma}}_{AB}$ when $m_A>m_B$
  and zero otherwise, we have
    $$
\sum_{k=1}^n
  X_k(N_k-X_k)(A_k-m_{A,k})(B_k-m_{B,k})/(n^{m_B/m_A}n^{1/2(1-m_B/m_A)})
  \stackrel{a.s.}\longrightarrow \widetilde{\Sigma}_{AB}.$$ Summing
  up, condition (c2) of Theorem \ref{thm:triangular} is satisfied.
  Regarding the last condition (c3), we observe that
  $|\mathbf{T}_{n,k}|=O(1/\sqrt{n})+O(1/\sqrt{n^{m_B/m_A}})
  \stackrel{a.s.}\longrightarrow
  0$. We are now ready to apply Theorem \ref{thm:triangular}, that
  gives the stable convergence of $\sum_{k=1}^n\mathbf{T}_{n,k}$
  toward $\mathcal{N}(0,\widetilde{\Sigma})$, where
  $$
\widetilde{\Sigma}=
  \begin{pmatrix}
    \widetilde{\Sigma}_{AA} &\widetilde{\Sigma}_{AB}\\
    \widetilde{\Sigma}_{AB} &\widetilde{\Sigma}_{BB}
  \end{pmatrix}.
  $$
In order to conclude, it is enough
  to observe that
  \begin{equation*}
    \begin{split}
      &  \sqrt{N_{A,n}}\left(\widehat{m}_{A,n}-m_A\right) =
  \sqrt{\frac{n}{N_{A,n}}}\frac{\sum_{k=1}^nX_kA_k}{\sqrt{n}}-\sqrt{N_{A,n}}m_A=\\
&  \sqrt{\frac{n}{N_{A,n}}}\frac{\sum_{k=1}^nX_k(A_k-m_{A,k})}{\sqrt{n}}
      +\sqrt{N_{A,n}}\left(\frac{\sum_{k=1}^n X_km_{A,k}}{N_{A,n}}-m_A\right)= \\
&  \sqrt{\frac{n}{N_{A,n}}}\sum_{k=1}^n [\mathbf{T}_{n,k}]_1
      +\sqrt{\frac{n}{N_{A,n}}}\frac{\sum_{k=1}^n X_k(m_{A,k}-m_A)}{n},
    \end{split}
  \end{equation*}
  and, similarly,
       \begin{equation*}
    \begin{split}
      & \sqrt{N_{B,n}}\left(\widehat{m}_{B,n}-m_B\right) =
      \sqrt{\frac{n^{m_B/m_A}}{N_{B,n}}}
      \frac{\sum_{k=1}^n(N_k-X_k)B_k}{\sqrt{n^{m_B/m_A}}}-\sqrt{N_{B,n}}m_B=\\ &
      \sqrt{\frac{n^{m_B/m_A}}{N_{B,n}}}
      \frac{\sum_{k=1}^n(N_k-X_k)(B_k-m_{B,k})}{\sqrt{n^{m_B/m_A}}}
      +\sqrt{N_{B,n}}\left(\frac{\sum_{k=1}^n
        (N_k-X_k)m_{B,k}}{N_{B,n}}-m_B\right)= \\ &
      \sqrt{\frac{n^{m_B/m_A}}{N_{B,n}}}\sum_{k=1}^n
           [\mathbf{T}_{n,k}]_2+
           \sqrt{\frac{n^{m_B/m_A}}{N_{A,n}}}\frac{\sum_{k=1}^n
             (N_k-X_k)(m_{B,k}-m_B)}{n^{m_B/m_A}}.
    \end{split}
  \end{equation*}
 Therefore, recalling that $N_{A,n}/n$ converges almost surely toward
 $N$ in the case $m_A>m_B$ and toward $NZ$ in the case $m_A=m_B$ and
 $N_{B,n}/n^{m_B/m_A}$ converges almost surely toward
 $\frac{m_A}{m_B}N\widetilde{Z}$ in the case $m_A>m_B$ and toward
 $N(1-Z)$ in the case $m_A=m_B$ and using Lemma \ref{lemma-app} in
 order to prove that $\sum_{k=1}^n X_k(m_{A,k}-m_A)/n$ and
 $\sum_{k=1}^n (N_k-X_k)(m_{B,k}-m_B)/n^{m_B/m_A}$ converges almost surely
 toward zero, we can conclude.
 \end{proof}

\begin{corollary}\label{normalization-th}
With the same notation and assumptions as in
Theorem~\ref{asymptotic-normality},  set $\sigma^2_{A}=q_{A}-m_{A}^2$,
$\sigma^2_{B}=q_{B}-m_{B}^2$ and
$\rho_{AB}=c_{AB}/\sqrt{\sigma_A^2\sigma_B^2}=
(q_{AB}-m_{A}m_{B})/\sqrt{\sigma^2_{A}\sigma_B^2}$.
Moreover, assume $\sigma_A^2>0$ and $\sigma_B^2>0$ and define
  $$
  \lambda_n=
  \frac{\sigma_A^2N_{B,n}/n}{\sigma^2_AN_{B,n}/n+\sigma_B^2N_{A,n}/n}
  $$
  and
  $$
  \Gamma_n=\begin{cases}
  1+\left(\frac{Q}{N}-1\right)\lambda_n\quad&\mbox{if }
  m_A>m_B\\
  \frac{Q}{N}+\left(\frac{Q}{N}-1\right)\left[(2Z-1)\lambda_n
    +2\rho_{AB}\sqrt{\lambda_n(1-\lambda_n)}-Z\right]\quad&\mbox{if } m_A=m_B.
  \end{cases}
$$
Then, in both cases, we have
  $$
  \zeta_n=
\frac{1}{\sqrt{\Gamma_n}}
\frac{\widehat{m}_{A,n}-\widehat{m}_{B,n}-(m_A-m_B)}
       { \sqrt{\sigma^2_{A,n}/N_{A,n}+\sigma^2_{B,n}/N_{B,n}} }
       \stackrel{stably}\longrightarrow\mathcal{N}(0,1)\,,
       $$
provided $\Gamma=a.s.-\lim_n\Gamma_n>0$ almost surely.
\end{corollary}

\begin{proof}
In both cases, we have  
  $$
  \lambda_n=
  \frac{\sigma_A^2N_{B,n}/n}{\sigma^2_AN_{B,n}/n+\sigma_B^2N_{A,n}/n}
  \stackrel{a.s.}\longrightarrow
  \lambda=\frac{\sigma_A^2(1-Z)}{\sigma^2_A(1-Z)+\sigma_B^2Z},
  $$ where $Z$ is equal to $1$ (and so $\lambda=0$) in the case
  $m_A>m_B$ and it is a non-atomic random variable with values in
  $[0,1]$ in the case $m_A=m_B$.  Moreover, we have
  $$
  \Gamma_n\stackrel{a.s.}\longrightarrow \Gamma=
\begin{cases}
  1+\left(\frac{Q}{N}-1\right)\lambda=1\quad&\mbox{if }
  m_A>m_B\\
  \frac{Q}{N}+\left(\frac{Q}{N}-1\right)\left[(2Z-1)\lambda
    +2\rho_{AB}\sqrt{\lambda(1-\lambda)}-Z\right]\quad&\mbox{if } m_A=m_B,
  \end{cases}
$$
that is
$$
\Gamma_n\stackrel{a.s.}\longrightarrow
\Gamma=\frac{\lambda}{\sigma_A^2}\Sigma_{AA}+
       \frac{1-\lambda}{\sigma_B^2}\Sigma_{BB}+
       2\frac{\sqrt{\lambda(1-\lambda)}}{\sqrt{\sigma^2_A\sigma^2_B}}\Sigma_{AB}.
  $$
  Therefore, in order to conclude it is enough to note that 
  $$
  \zeta_n=\frac{1}{\sqrt{\Gamma_n}}
  \left(
  \sqrt{\lambda_n}\frac{\sqrt{N_{A,n}}}{\sqrt{\sigma^2_A}}
  (\widehat{m}_{A,n}-m_A)+
  \sqrt{1-\lambda_n}\frac{\sqrt{N_{B,n}}}{\sqrt{\sigma^2_B}}
  (\widehat{m}_{B,n}-m_B)
  \right).
  $$
  Indeed, by the above Theorem \ref{asymptotic-normality} and the
  above almost sure convergences, we get that
$$
\zeta_n\stackrel{stably}\longrightarrow
\mathcal{N}(0,1).
$$
\end{proof}

\begin{remark}\rm 
  Note that, in the case $m_A=m_B$, since $Z$ is a non-atomic random
  variable with values in $[0,1]$, the limit random variable $\Gamma$
  is equal to $0$ with a strictly positive probability if and only if
  $P(Q>N)=P(N>1)>0$ and, on the event $\{N>1\}$, the random variable
  $Q/N$ is a function of $Z$ such that
  $$
-\frac{Q}{N}/\left(\frac{Q}{N}-1\right)=\left[(2Z-1)\lambda
    +2\rho_{AB}\sqrt{\lambda(1-\lambda)}-Z\right],
$$
hence, in a very special case.
\end{remark}


\section{test}\label{test}
Set $M_n=\frac{1}{n}\sum_{j=1}^nX_j/N_j$, which is the empirical mean
of the proportions of balls of color $A$ in the samples until
time-step $n$. By Remark 3.6 in \cite{cri-lou-min-multidrawing}, we
know that $M_n\stackrel{a.s.}\longrightarrow Z$. Moreover, set
$\widehat{\sigma}^2_{A,n}=\widehat{q}_{A,n}-\widehat{m}_{A,n}^2$,
$\widehat{\sigma}^2_{B,n}=\widehat{q}_{B,n}-\widehat{m}_{B,n}^2$ and
$\widehat{\rho}_{AB,n}=
(\widehat{q}_{AB,n}-\widehat{m}_{A,n}\widehat{m}_{B,n})/
\sqrt{\widehat{\sigma}^2_{A,n}\widehat{\sigma}^2_{B,n}}$.  The stable
convergence stated in Corollary \ref{normalization-th} still holds
true even if we replace all the quantities with their strongly
consistent estimators. More precisely, in both cases, assuming
$\sigma_A^2>0$, $\sigma_B^2>0$ and $\Gamma=a.s.-\lim_n \Gamma_n>0$
almost surely, we have
$$
  \widehat{\zeta}_n=
\frac{1}{\sqrt{\widehat{\Gamma}_n}}
\frac{\widehat{m}_{A,n}-\widehat{m}_{B,n}-(m_A-m_B)}
       { \sqrt{\widehat{\sigma}^2_{A,n}/N_{A,n}+\widehat{\sigma}^2_{B,n}/N_{B,n}} }
       \stackrel{stably}\longrightarrow\mathcal{N}(0,1)\,,
       $$ where $\widehat{\Gamma}_n$ is defined as $\Gamma_n$, but
       replacing the the random variables $Z$, $N$ and $Q$ by $M_n$,
       $\widehat{\mu}_n$ and $\widehat{q}_{N,n}$, respectively, and
       replacing the quantities $\sigma^2_A,\,\sigma^2_B$ and
       $\rho_{AB}$ by their estimators
       $\widehat{\sigma}^2_{A,n},\,\widehat{\sigma}^2_{B,n}$ and
       $\widehat{\rho}_{AB,n}$. Given this fact, a
       critical region (with asymptotic level $\theta$) for the
       hypothesis test on the reinforcement means
$$
H_0: m_A = m_B\qquad\mbox{versus}\qquad
H_1: m_A > m_B
$$
is given by 
$$
C_\theta= \{\widehat{\zeta}^0_n> q_{1-\theta}\},
$$ where $\widehat{\zeta}^0_n$ is defined as the random variable
$\widehat{\zeta}_n$ in the case $m_A=m_B$ and $q_{1-\theta}$ is the
quantile of the standard normal distribution of order $1-\theta$,
that is $\mathcal{N}(0,1)(q_{1-\theta},+\infty)=\theta$.

\begin{remark}(Power of the test)\\
  \rm 
  We have 
  \begin{equation*}
    \begin{split}
  \widehat{\zeta}^0_n &=
  \frac{1}{\sqrt{\widehat{\Gamma}^{0}_n}}
  \frac{\widehat{m}_{A,n}-\widehat{m}_{B,n}}
       { \sqrt{ \widehat{\sigma}^2_{A,n}/N_{A,n}+\widehat{\sigma}^2_{B,n}/N_{B,n}} }
       \\
       &=
       \sqrt{\frac{\widehat{\Gamma}^{\neq 0}_n}{\widehat{\Gamma}^0_n}}
       \frac{1}{\sqrt{\widehat{\Gamma}^{\neq 0}_n}}
  \frac{\widehat{m}_{A,n}-\widehat{m}_{B,n}-(m_A-m_B)}
      { \sqrt{ \widehat{\sigma}^2_{A,n}/N_{A,n}+\widehat{\sigma}^2_{B,n}/N_{B,n}}}
       +
       \frac{1}{\sqrt{\widehat{\Gamma}^0_n}}
       \frac{m_A-m_B}
       { \sqrt{ \widehat{\sigma}^2_{A,n}/N_{A,n}+\widehat{\sigma}^2_{B,n}/N_{B,n}}}
       \,,
    \end{split}
    \end{equation*}
    where $\widehat{\Gamma}_n^0$ (resp. $\widehat{\Gamma}_n^{\neq 0}$)
    is the random variable defined as $\Gamma_n$ in the case $m_A=m_B$
    (resp. $m_A>m_B$), but replacing the random variables $Z$, $N$ and
    $Q$ by $M_n$, $\widehat{\mu}_n$ and $\widehat{q}_{N,n}$,
    respectively, and the quantities $\sigma^2_A,\,\sigma^2_B$ and
    $\rho_{AB}$ by their estimators
    $\widehat{\sigma}^2_{A,n},\,\widehat{\sigma}^2_{B,n}$ and
    $\widehat{\rho}_{AB,n}$. Under the alternative hypothesis
    (i.e. $m_A>m_B$), we have $Z=1$ and so $\lambda=0$ and,
    consequently, $\widehat{\Gamma}^0_n\stackrel{a.s.}\longrightarrow
    1$ and $\widehat{\Gamma}^{\neq 0}_n\stackrel{a.s.}\longrightarrow
    1$ and, hence, by Corollary~\ref{normalization-th}, we get
    $$
    \sqrt{\frac{\widehat{\Gamma}^{\neq 0}_n}{\widehat{\Gamma}^0_n}}
       \frac{1}{\sqrt{\widehat{\Gamma}^{\neq 0}_n}}
  \frac{\widehat{m}_{A,n}-\widehat{m}_{B,n}-(m_A-m_B)}
      { \sqrt{ \widehat{\sigma}^2_{A,n}/N_{A,n}+\widehat{\sigma}^2_{B,n}/N_{B,n}}}
\stackrel{stably}\longrightarrow \mathcal{N}(0,1).
$$
Therefore, we can say that, under the alternative hypothesis, we have
  $$
  \widehat{\zeta}^0_n\stackrel{d}\approx
  \mathcal{N}\left(
\frac{1}{\sqrt{\widehat{\Gamma}^0_n}}
       \frac{m_A-m_B}
     { \sqrt{ \widehat{\sigma}^2_{A,n}/N_{A,n}+\widehat{\sigma}^2_{B,n}/N_{B,n}}}, 
  1\right).
  $$
  Hence, the power of the test can be approximated by
  \begin{equation}\label{approxpowerquantilenor}
  \mathcal{N}\left(
\frac{1}{\sqrt{\widehat{\Gamma}^0_n}}
\frac{\widehat{m}_{A,n}-\widehat{m}_{B,n}}
     { \sqrt{ \widehat{\sigma}^2_{A,n}/N_{A,n}+\widehat{\sigma}^2_{B,n}/N_{B,n}}},
     1\right)(q_{1-\theta},+\infty)\,.
     \end{equation}
Note that, under the alternative hypothesis, by Theorem~\ref{rate-th}, we have
$$
n^{m_B/m_A}
\left( \widehat{\sigma}^2_{A,n}/N_{A,n}+\widehat{\sigma}^2_{B,n}/N_{B,n}\right)
  \stackrel{a.s.}\longrightarrow
  \sigma^2_B \frac{m_B}{m_A}(N\widetilde{Z})^{-1}
  $$ and so the average value in the above normal distribution
  goes to $+\infty$ as $\sqrt{n^{m_B/m_A}}$.
\end{remark}

\begin{remark}(A very special case)\\
  \rm Note that the very special case when $N_n=1$ and $m_{A,n}=m_A$
  and $m_{B,n}=m_B$ for each $n$ corresponds to the setting studied in
  \cite{mayflournoy}.
\end{remark}

\begin{remark}(Some possible simplifications)\\
  \rm When $N$ and $Q$ are known, we do not need to use their
  estimators \eqref{def-est-NQ}. For example, a case is when all the
  sample sizes $N_n$ are equal to the same known
  constant~$\kappa$. More generally, it might be the case that the
  generating mechanism of the sample sizes $N_n$ is decided by the
  experimenter (for example, this is the case of a response-adaptive
  design) and so it might be that the random variables $N$ and $Q$ can
  be written explicitly.
\end{remark}

\begin{figure}[h!]
\centering 
\includegraphics[scale=0.6,keepaspectratio=true]{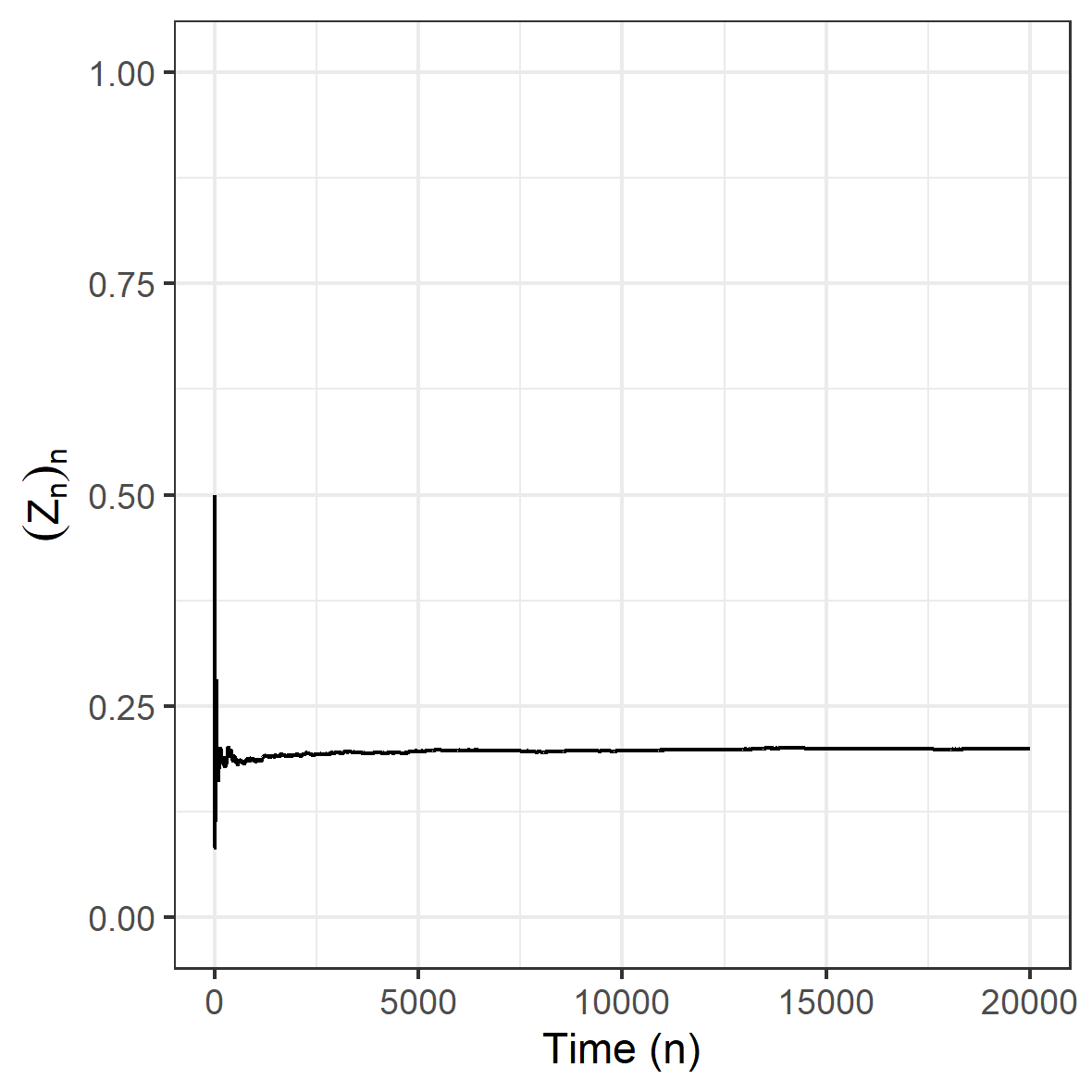}
\includegraphics[scale=0.6,keepaspectratio=true]{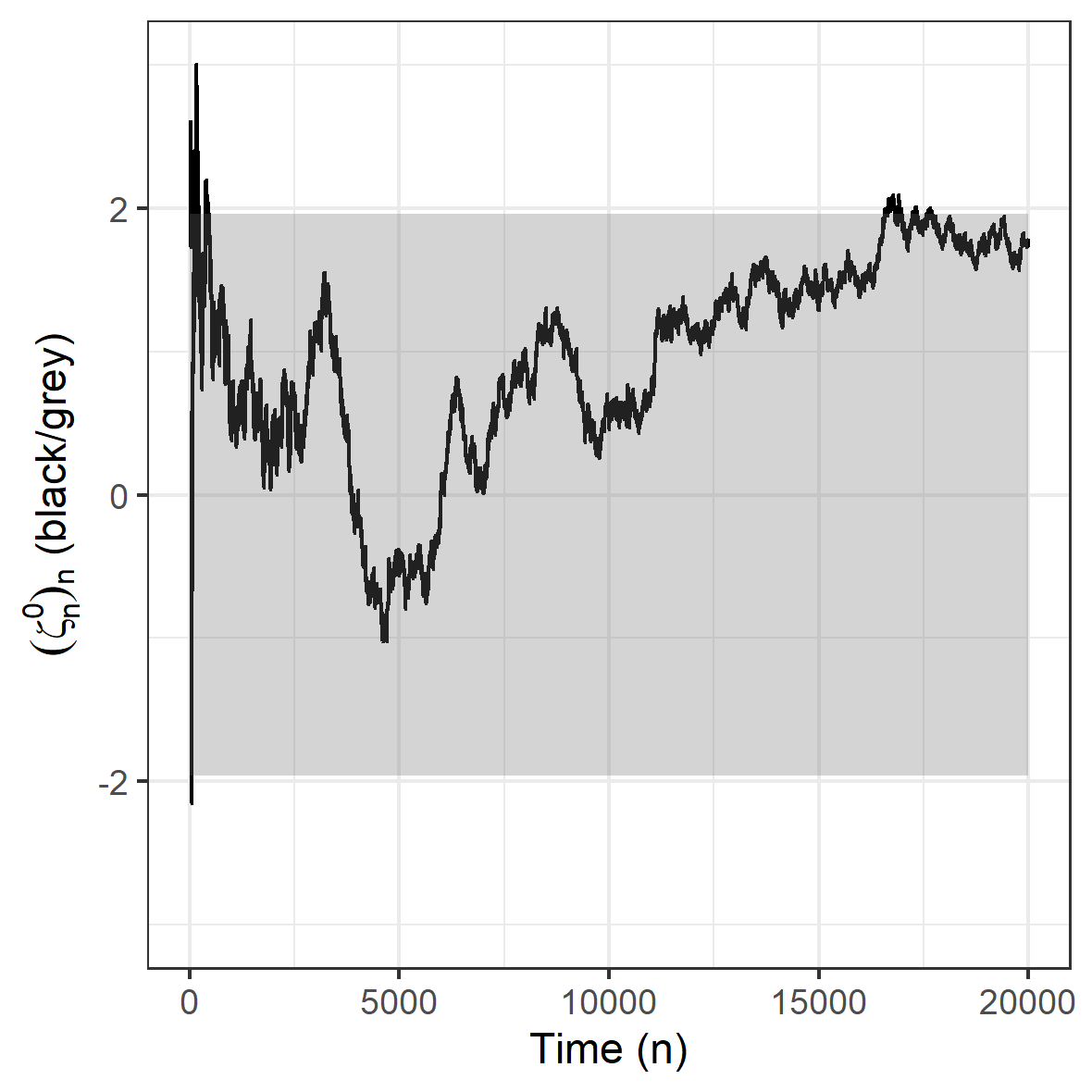}
\caption{Multinomial case under $H_0$ ($m_A=m_B$), with size$=12$,
  $p_A=4/15$, $p_B=4/15$. One sample path for $(Z_n)_n$ and
  $(\zeta^0_n)_n$. The gray region corresponds to a confidence
  interval of level $5\%$ given by the standard normal distribution.}
\label{fig:caso1b:H0}
\end{figure}

\noindent \textbf{Example}\\
Take each $N_n$ independent of $\mathcal{F}_{n-1}$ and uniformly
distributed on $\{1,\hdots,5\}$.  Moreover, take $[A_n,B_n]$ such that
 $$
A_n\stackrel{d}=1+Y_1\qquad\mbox{and}\qquad B_n\stackrel{d}=1+Y_2\,,
$$ where $Y_1$ and $Y_2$ are, respectively, the first and the second
component of a multinomial distribution associated to the parameters:
size$=12$, probabilities$=(p_A,p_B,p_3)$ with $p_A+p_B+p_3=1$. Thus
the random variables $A_n$ and $B_n$ are negatively correlated. We set
$a=b=5$.  \\ \indent Choosing $p_A=p_B=4/15$, the null hypothesis
$m_A=m_B$ holds. Fig.~\ref{fig:caso1b:H0} shows a sample path of
$(Z_n)_n$ and of $(\zeta^0_n)_n$. The number of iterations is
$n=20\,000$.

\begin{figure}[h!]
\centering 
\includegraphics[scale=0.8,keepaspectratio=true]{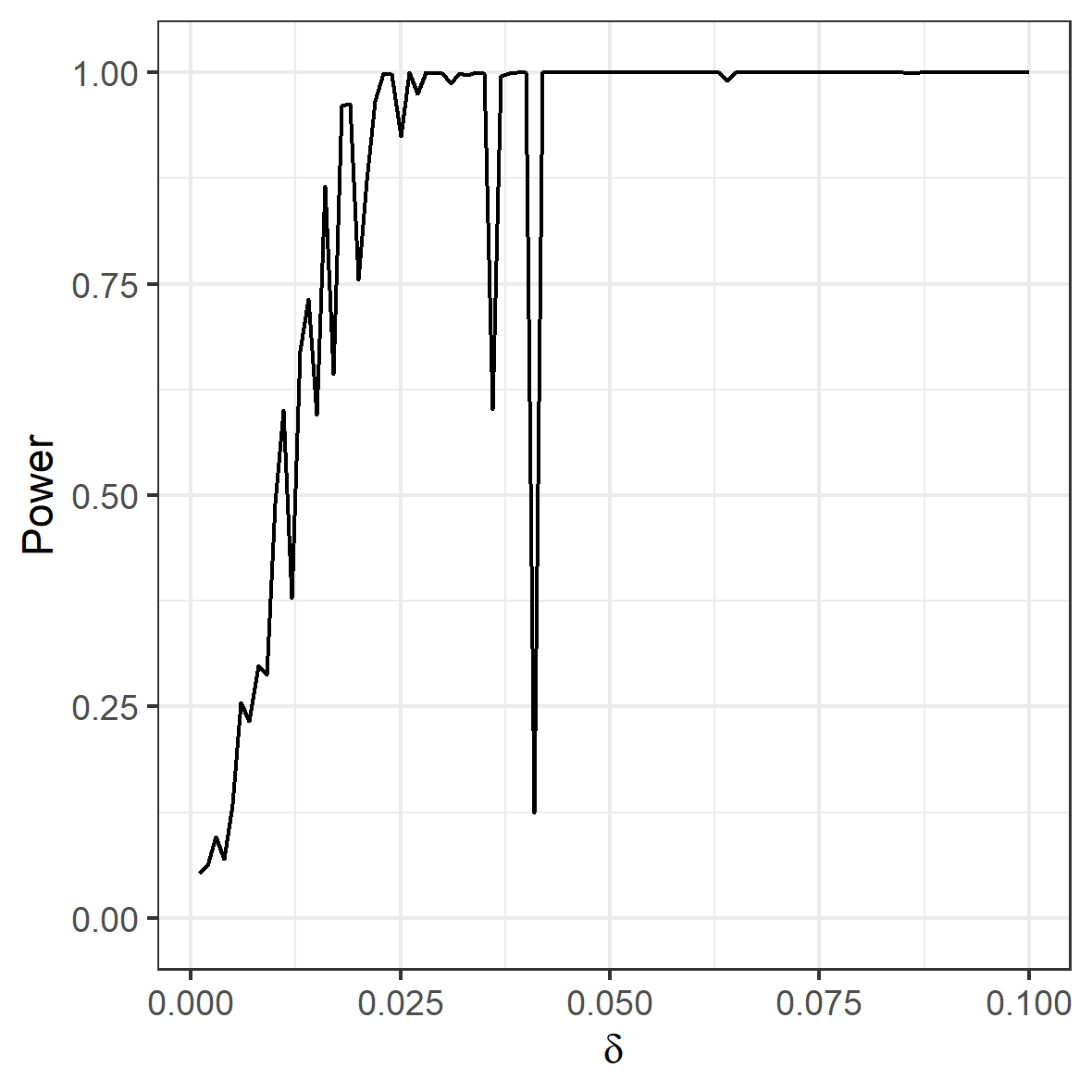}
\caption{Approximated power of the test in the Multinomial case, with
  size $=12$ and $p_B=4/15$. The parameters $p_A$ and $p_3$ vary and
  $\delta=m_A-m_B>0$. }
\label{fig:caso1b:H1:power}
\end{figure}

In order to obtain an empirical power, we make $\delta=m_A-m_B>0$ vary
(from $0$ to $0.1$ by $0.005$). We keep $p_B=4/15$ and the size$=12$
and we vary $p_A$ and $p_3$ for getting the different values of
$\delta$. We run the dynamics until time-step $n=10\,000$ and we use
the approximated power given by~\eqref{approxpowerquantilenor}.  Results
are sum up in Fig.~\ref{fig:caso1b:H1:power}.

\appendix

\section{Technical results}\label{app-technical}

\begin{lemma}\label{supermart}
    If $m_A>m_B$ and condition \eqref{cond-N} is satisfied, then, we
    have $1/K_n=o(n^{-1/\theta})$ for each $\theta>m_A/m_B$.
\end{lemma}
\begin{proof} Take $\theta>m_A/m_B$. We observe that
  \begin{equation*}
    \begin{split}
      &E\left[\frac{H_{n+1}}{K_{n+1}^\theta}-
        \frac{H_n}{K_n^\theta}\,|\,\mathcal{H}_n\right]=
      E\left[\frac{H_{n+1}}{K_n^\theta}-\frac{H_n}{K_n^\theta} +
        \frac{H_{n+1}}{K_{n+1}^\theta}-\frac{H_{n+1}}{K_n^\theta}\,|\,\mathcal{H}_n
        \right]=\\
      &\sum_{k\in\mathcal{X}_{n+1}}
      p_{n+1,k}\left(\frac{H_n+A_{n+1}k}{K_n^\theta}-\frac{H_n}{K_n^\theta}\right)
      +\sum_{k\in\mathcal{X}_{n+1}}
      p_{n+1,k}(H_n+A_{n+1}k)
      \left(\frac{1}{(K_n+B_{n+1}(N_{n+1}-k))^\theta}-\frac{1}{K_n^\theta}\right)=\\
      &
\sum_{k\in \mathcal{X}_{n+1}} p_{n+1,k}\frac{A_{n+1}k}{K_n^\theta}
+
\sum_{k\in \mathcal{X}_{n+1}\setminus \{N_{n+1}\}} p_{n+1,k}(H_n+A_{n+1}k)
      \left(\frac{1}{(K_n+B_{n+1}(N_{n+1}-k))^\theta}-\frac{1}{K_n^\theta}\right)
      \,.
      \end{split}
  \end{equation*}
  Using the Taylor expansion of the function $f(x)=1/(c+x)^\theta$ with
  $c=K_n$ and $x=B_{n+1}(N_{n+1}-k)$, we can choose a constant
  $\xi$ such that eventually
  $$
  \left(\frac{1}{(K_n+B_{n+1}(N_{n+1}-k))^{\theta}}-\frac{1}{K_n^\theta}\right)
  \leq
  -\frac{\theta}{K_n^{\theta+1}}\left(B_{n+1}(N_{n+1}-k)-\frac{\xi}{K_n}\right)\,.
  $$ Therefore the last term of the above equalities is eventually
  smaller than or equal to
  $$
  \frac{H_n}{K_n^\theta}\left\{
  \sum_{k\in \mathcal{X}_{n+1}}
 \left( \frac{A_{n+1}k}{H_n}- \theta\frac{B_{n+1}(N_{n+1}-k)}{K_n} \right)p_{n+1,k}
 + \theta\xi \sum_{k\in \mathcal{X}_{n+1}\setminus\{N_{n+1}\}}
 \frac{(1+A_{n+1}k/H_n)}{K_n^2}p_{n+1,k}\right\}\,.
$$
Now, we observe that we have 
\begin{equation*}
  \begin{split}
  E\left[\sum_{k\in \mathcal{X}_{n+1}}
    \left(\frac{A_{n+1}k}{H_n}-\theta \frac{B_{n+1}(N_{n+1}-k)}{K_n}\right)p_{n+1,k}
    \,|\,\mathcal{G}_n\right]
  =
  \frac{N_{n+1}}{S_n}(m_{n+1,A}-m_{n+1,B}\theta)
  \end{split}
  \end{equation*}
and, by Lemma~A.1 in \cite{cri-lou-min-multidrawing}, 
$$
E\left[\sum_{k\in \mathcal{X}_{n+1}\setminus\{N_{n+1}\}}
  \frac{(1+A_{n+1}k/H_n)}{K_n^2}p_{n+1,k}\,|\,\mathcal{G}_n\right]
\leq \frac{(1-p_{n+1,N_{n+1}})}{K_n^2}+\frac{m_{n+1,A}N_{n+1}}{S_n K_n^2}
=O(1/(S_n K_n))\,.
$$
Therefore, since $N_{n+1}\geq 1$, we have 
$$
E\left[\frac{H_{n+1}}{K_{n+1}^\theta}-\frac{H_n}{K_n^\theta}\,|\,\mathcal{G}_n\right]
\leq
\frac{H_n}{K_n^\theta}\frac{N_{n+1}}{S_n}
\left[-(m_{n+1,B}\theta-m_{n+1,A})+O(1/K_n)\right]
$$ and so, for $\theta>m_A/m_B$, since $K_n\uparrow +\infty$ (by Lemma
3.1 in \cite{cri-lou-min-multidrawing}), we can conclude that the
above conditional expectation is eventually negative. This fact means
that $H_n/K_n^\theta$ is eventually a positive supermartingale and so
it converges almost surely toward a finite random variable. Since
$H_n/S_n$ converges toward $1$ almost surely and $S_n$ converges
toward $m_AN>0$ almost surely, we have that $H_n/n$ converges almost
surely toward a finite random variable and so also $n/K_n^{\theta}$
converges almost surely toward a finite random variable for each
$\theta>m_A/m_B$. This means that $n/K_n^{\theta}$ converges almost
surely toward $0$ for each $\theta>m_A/m_B$, that is
$1/K_n=o(n^{1/\theta})$ for each $\theta>m_A/m_B$.
\end{proof}
  
  \begin{lemma}\label{lemma-limit}
If we are in case $m_A>m_B$ and condition \eqref{cond-N} is
satisfied, then $(H_n/K_n^{m_A/m_B})_n$ converges almost surely toward
a random variable taking values in $(0,+\infty)$.
\end{lemma}

 \begin{proof} 
  Set $L_n=\ln(H_n/K_n^{m_A/m_B})$.  If we prove that $L_n$ converges
  almost surely to a finite random variable (see Lemma
  \ref{lemma-pemantle}), then $Y_n=H_n/K_n^{m_A/m_B}$ converges to a
  random variable $Y$ with values in $(0,+\infty)$. In order to prove
  the almost sure convergence of $(L_n)$, we are going to apply Lemma
  \ref{lemma-pemantle}. Therefore, we set
  $\Delta_n=E[L_{n+1}-L_n|\mathcal{G}_n]$ and
  $Q_n=E[(L_{n+1}-L_n)^2|\mathcal{G}_n]$. We recall that, from Lemma
  \ref{supermart} we know that $1/K_n=O(1/n^\gamma)$ for some $\gamma>0$
  and, moreover, we know that $H_n/n\stackrel{a.s.}1/(m_AN)>0$ and so
  $1/H_n=O(1/n)$. Using the notation \eqref{hypergeometric-distr}, we
  have
  \begin{equation*}
  \begin{split}
    &\Delta_n=E[\ln(H_{n+1})-\ln(H_n)|\mathcal{H}_n]-
    \frac{m_A}{m_B}E[\ln(K_{n+1})-\ln(K_n)|\mathcal{H}_n]=
    \\ &\sum_{k\in\mathcal{X}_{n+1}}
    \left\{\left(\ln(H_{n}+A_{n+1}k)-\ln(H_n)\right)- \frac{m_A}{m_B}
    \left(\ln(K_{n}+B_{n+1}(N_{n+1}-k))-\ln(K_n)\right)\right\}p_{n+1,k}=
    \\ &\sum_{k\in\mathcal{X}_{n+1}}
    \left\{\int_0^{A_{n+1}k}\frac{1}{H_n+t}\,dt-
    \frac{m_A}{m_B}\int_0^{B_{n+1}(N_{n+1}-k)}\frac{1}{K_n+t}\,dt\right\}p_{n+1,k}\,.
  \end{split}
  \end{equation*}
  Since $1/(H_n+t)\leq 1/H_n$ and $1/(K_n+t)\geq 1/K_n-t/K_n^2$ for
  each $t\geq 0$ and each $n$, the last term of the above equalities
  is eventually smaller than or equal to 
  $$
\sum_{k\in\mathcal{X}_{n+1}}\left\{\frac{A_{n+1}k}{H_n}-
\frac{m_A}{m_B}\frac{B_{n+1}(N_{n+1}-k)}{K_n}+
\frac{m_A}{m_B}\frac{B_{n+1}^2(N_{n+1}-k)^2}{2 K_n^2}\right\}p_{n+1,k}\,.
$$
Now, we observe that
\begin{equation*}
  \begin{split}
 & E[\sum_{k\in\mathcal{X}_{n+1}}\left(\frac{A_{n+1}k}{H_n}-
  \frac{m_A}{m_B}
  \frac{B_{n+1}(N_{n+1}-k)}{K_n}\right)p_{n+1,k}\,|\,\mathcal{G}_n]=
\frac{m_{n+1,A}N_{n+1}H_n}{H_nS_n}-\frac{m_{A}m_{B,n+1}}{m_B}\frac{N_{n+1}K_n}{K_nS_n}
=\\
&\frac{N_{n+1}}{S_n}\left(m_{A,n+1}-\frac{m_A}{m_B}m_{B,n+1}\right)
=O(|m_{A,n+1}m_B-m_A m_{B,n+1}|/n)\,.
  \end{split}
\end{equation*}
  Therefore, we have
  \begin{equation*}
    \begin{split}
\Delta_n&\leq O(|m_{A,n+1}m_B-m_A m_{B,n+1}|/n)+
\frac{C^3}{2 K_n^2}\sum_{k\in\mathcal{X}_{n+1}}
(N_{n+1}-k)p_{n+1,k}\\
&=
O(|m_{A,n+1}m_B-m_A m_{B,n+1}|/n)+
\frac{C^3}{2 K_n^2}N_{n+1}\frac{K_n}{S_n}
=
O(|m_{A,n+1}m_B-m_A m_{B,n+1}|/n) + O(1/(K_nS_n))\\
&=
O(|m_{A,n+1}m_B-m_A m_{B,n+1}|/n)+
O(1/n^{1+\gamma})\,.
    \end{split}
\end{equation*}
    Finally, we note that $-\Delta_n=\ln(K_{n+1}^{m_B/m_A}/H_{n+1})-
\ln(K_n^{m_B/m_A}/H_n)$
and so, with the same arguments as before, we get 
$-\Delta_n\leq O(|m_{A,n+1}m_B-m_A m_{B,n+1}|/n)+O(1/n^{1+\gamma})$.
Thus, $\sum_n|\Delta_n|<+\infty$ almost surely.  Similarly, we have
 \begin{equation*}
  \begin{split}
&E[(\ln(H_{n+1})-\ln(H_n)-\frac{m_A}{m_B}\ln(K_{n+1})+\ln(K_n))^2|\mathcal{H}_n]
\leq\\
&2\left\{
E[(\ln(H_{n+1})-\ln(H_n))^2|\mathcal{H}_n]+
\frac{m_A^2}{m_B^2}E[(\ln(K_{n+1})-\ln(K_n))^2|\mathcal{H}_n]
\right\}
\leq\\
&2
\sum_{k\in\mathcal{X}_{n+1}}
\left(\frac{A_{n+1}^2k^2}{H_n^2}+
\frac{m_A^2}{m_B^2}\frac{B_{n+1}^2(N_{n+1}-k)^2}{K_n^2}\right)p_{n+1,k}
\leq
\\
&2C^3
\sum_{k\in\mathcal{X}_{n+1}}
\left(\frac{k}{H_n^2}+\frac{m_A^2}{m_B^2}\frac{(N_{n+1}-k)}{K_n^2}\right)p_{n+1,k}
\leq 2C^4\left(\frac{1}{H_nS_n}+\frac{m_A^2}{m_B^2}\frac{1}{K_nS_n}\right)\,.
\end{split}
 \end{equation*}
 Therefore, we get $Q_n=O(1/n^{1+\gamma})$ and so $\sum_nQ_n<+\infty$
 almost surely.
  \end{proof}
  

 \section{Some auxiliary results}\label{app-auxiliary}

We recall two known technical results:

\begin{lemma}\label{lemma-app} (Lemma 4.1 in \cite{cri-dai-min})\\
Let $(\alpha_j)$ and $(\beta_j)_j$ be two sequences of strictly
positive numbers with
$$\beta_n\uparrow +\infty\qquad\mbox{ and }\qquad
\frac{1}{\beta_n}\sum_{j=1}^n\frac{1}{\alpha_j}\longrightarrow \gamma.
$$ Let $(Y_n)$ be a
sequence of real random variables, adapted to a filtration $\mathcal
F$. If $\sum_{j\geq 1} E[Y_j^2]/(\alpha_j\beta_j)^2<+\infty$ and
$E[Y_{j}| {\mathcal F}_{j-1}]\stackrel{a.s.}\longrightarrow Y$ for
some real random variable $Y$, then
$$
\frac{1}{\beta_n}\sum_{j=1}^n \frac{Y_j}{\alpha_j}\stackrel{a.s.}\longrightarrow
\gamma Y.
$$
\end{lemma}

\begin{lemma}\label{lemma-pemantle} (Lemma~3.2 in~\cite{pemantle-volkov-1999})\\
  Let $(L_n)$ be a sequence of random variables, adapted to a
  filtration $\mathcal{G}_n$. Set
  $\Delta_n=E[L_{n+1}-L_n|\mathcal{G}_n]$ and
  $Q_n=E[(L_{n+1}-L_n)^2|\mathcal{G}_n]$. If $\sum_n\Delta_n$ and
  $\sum_nQ_n$ are almost surely convergent, then $(L_n)$ converges
  almost surely to a finite random variable.
\end{lemma}


\section{Stable convergence}\label{app-stable-conv}
This brief appendix contains some basic definitions and results
concerning stable convergence. For more details, we
refer the reader to \cite{crimaldi-libro, hall-1980} and the
references therein.\\

\indent Let $(\Omega, {\mathcal A}, P)$ be a probability space,
and let $S$ be a Polish space, endowed with its Borel
$\sigma$-field. A {\em
  kernel} on $S$, or a random probability measure on $S$, is a
collection $K=\{K(\omega):\, \omega\in\Omega\}$ of probability
measures on the Borel $\sigma$-field of $S$ such that, for each
bounded Borel real function $f$ on $S$, the map
$$
\omega\mapsto K\!f(\omega)=\int f (x)\, K(\omega)(dx)
$$
is $\mathcal A$-measurable. Given a sub-$\sigma$-field $\mathcal
H$ of $\mathcal A$, a kernel $K$ is said $\mathcal H$-measurable
if all the
above random variables $K\!f$ are $\mathcal H$-measurable.\\

\indent On $(\Omega, {\mathcal A},P)$, let $(Y_n)_n$ be a sequence
of $S$-valued random variables, let $\mathcal H$ be a
sub-$\sigma$-field of $\mathcal A$, and let $K$ be a $\mathcal
H$-measurable kernel on $S$. Then we say that $Y_n$ converges {\em
$\mathcal H$-stably} to $K$, and we write $Y_n\longrightarrow K$
${\mathcal H}$-stably, if
$$
P(Y_n \in \cdot \,|\, H)\stackrel{weakly}\longrightarrow
E\left[K(\cdot)\,|\, H \right] \qquad\hbox{for all } H\in{\mathcal
H}\; \hbox{with } P(H) > 0,
$$where $K(\cdot)$ denotes the random variable
  defined, for each Borel set $B$ of $S$, as $\omega\mapsto
  K\!I_B(\omega)=K(\omega)(B)$.  In the case when ${\mathcal
  H}={\mathcal A}$, we simply say that $Y_n$ converges {\em stably} to
$K$ and we write $Y_n\longrightarrow K$ stably. Clearly, if
$Y_n\longrightarrow K$ ${\mathcal H}$-stably, then $Y_n$ converges
in distribution to the probability distribution $E[K(\cdot)]$.
Moreover, the $\mathcal H$-stable convergence of $Y_n$ to $K$ can
be stated in terms of the following convergence of conditional
expectations:
\begin{equation}\label{def-stable}
E[f(Y_n)\,|\, {\mathcal H}]\stackrel{\sigma(L^1,\,
L^{\infty})}\longrightarrow K\!f
\end{equation}
for each bounded continuous real function $f$ on $S$. \\

\indent We denote by $\mathcal{N}(\mathbf{\mu},C)$ the multivariate
Gaussian probability distribution with vector of mean values
$\mathbf{\mu}$ and covariance matrix $C$. Therefore, when $\Sigma$ is
a random positive semidefinite matrix, the symbol
$\mathcal{N}(0,\Sigma)$ denotes the Gaussian kernel
$\{\mathcal{N}(0,\Sigma(\omega)):\, \omega\in\Omega\}\}$.
\\

From Proposition~3.1 in \cite{cri-pra-multivariate-mart}, we
can get the following result.

\begin{theorem}\label{thm:triangular}
Let $({\mathbf T}_{n,k})_{n\geq 1, 1\leq k\leq k_n}$ be a triangular
array of $d$-dimensional real random vectors, such that, for each
fixed $n$, the finite sequence $({\mathbf T}_{n,k})_{1\leq k\leq k_n}$
is a martingale difference array with respect to a given filtration
$({\mathcal G}_{n,k})_{k\geq 0}$. Moreover, let $(t_n)_n$ be a
sequence of real numbers and assume that the following
conditions hold:
\begin{itemize}
\item[(c1)] ${\mathcal G}_{n,k}{\underline{\subset}} {\mathcal G}_{n+1,
  k}$ for each $n$ and $1\leq k\leq k_n$;
\item[(c2)] $\sum_{k=1}^{k_n} (t_n{\mathbf
  T}_{n,k})(t_n{\mathbf T}_{n,k})^{\top}=t_n^2\sum_{k=1}^{k_n} {\mathbf
  T}_{n,k}{\mathbf T}_{n,k}^{\top} \stackrel{P}\longrightarrow \Sigma$,
  where $\Sigma$ is a random positive semidefinite matrix;
\item[(c3)] $\sup_{1\leq k\leq k_n} |t_n{\mathbf T}_{n,k}|
\stackrel{L^1}\longrightarrow 0$.
\end{itemize}
Then $t_n\sum_{k=1}^{k_n}{\mathbf T}_{n,k}$ converges stably to the
Gaussian kernel ${\mathcal N}(\mathbf{0}, \Sigma)$.
\end{theorem}
\bigskip


\noindent {\bf Acknowledgments}\\
\noindent Irene Crimaldi and Ida Minelli are members of the Italian
Group ``Gruppo Nazionale per l'Analisi Matematica, la Probabilit\`a e
le loro Applicazioni'' of the Italian Institute ``Istituto Nazionale
di Alta Matematica''.
\medskip 

\noindent {\bf Funding Sources}\\
\noindent Irene Crimaldi is partially supported by the Italian
``Programma di Attivit\`a Integrata'' (PAI), project ``TOol for
Fighting FakEs'' (TOFFE) funded by IMT School for Advanced Studies
Lucca.
\medskip 

\noindent {\bf Declaration}\\
\noindent All the authors equally contributed to this work.


\begin{thebibliography}{10}

\bibitem{AguechLasmarSelmi2019}
R.~Aguech, N.~Lasmar, and O.~Selmi.
\newblock A generalized urn with multiple drawing and random addition.
\newblock {\em Annals of the Institute of Statistical Mathematics},
  71(2):389--408, Apr. 2019.

\bibitem{AguechSelmi-unbalanced}
R.~Aguech and O.~Selmi.
\newblock Unbalanced multi-drawing urn with random addition matrix.
\newblock {\em Arab Journal of Mathematical Sciences}, 2019.

\bibitem{perron}
D.~A. Aoudia and F.~Perron.
\newblock A new randomized {P}\'olya urn model.
\newblock {\em Appl. Math.}, 3:2118--2122, 2012.

\bibitem{Chen2020}
M.-R. Chen.
\newblock A time dependent {P}{\'{o}}lya urn with multiple drawings.
\newblock {\em Probab. Eng. Informational Sci.}, 34(4):469--483, 2020.

\bibitem{chen-kuba}
M.-R. Chen and M.~Kuba.
\newblock On generalized {P}\'olya urn models.
\newblock {\em J. Appl. Prob.}, 50:1169--1186, 2013.

\bibitem{chen-wei-2005}
M.-R. Chen and C.-Z. Wei.
\newblock A new urn model.
\newblock {\em Journal of Applied Probability}, 42(4):964--976, Dec. 2005.

\bibitem{cri-ipergeom}
I.~Crimaldi.
\newblock Central limit theorems for a hypergeometric randomly reinforced urn.
\newblock {\em J. Appl. Prob.}, 53(3):899--913, 2016.

\bibitem{crimaldi2016}
I.~Crimaldi.
\newblock Central limit theorems for a hypergeometric randomly reinforced urn.
\newblock {\em Journal of Applied Probability}, 53(3):899--913, 2016.

\bibitem{crimaldi-libro}
I.~Crimaldi.
\newblock {\em {I}ntroduzione alla nozione di convergenza stabile e sue
  varianti ({I}ntroduction to the notion of stable convergence and its
  variants)}, volume~57.
\newblock {U}nione {M}atematica {I}taliana, Monograf s.r.l., Bologna, Italy.,
  2016.
\newblock Book written in Italian.

\bibitem{cri-dai-min}
I.~Crimaldi, P.~Dai~Pra, and I.~G. Minelli.
\newblock Fluctuation theorems for synchronization of interacting {P}\'olya's
  urns.
\newblock {\em Stochastic Process. Appl.}, 126(3):930--947, 2016.

\bibitem{cri-lou-min-multidrawing}
I.~Crimaldi, P.-Y. Louis, and I.~G. Minelli.
\newblock An urn model with random multiple drawing and random addition.
\newblock {\em Stochastic Processes and their Applications}, 147(C):270--299,
  2022.

\bibitem{cri-pra-multivariate-mart}
I.~Crimaldi and L.~Pratelli.
\newblock Convergence results for multivariate martingales.
\newblock {\em Stochastic Processes and their Applications}, 115(4):571--577,
  2005.

\bibitem{egg-pol}
F.~Eggenberger and G.~P\'{o}lya.
\newblock {\"{U}}ber die {S}tatistik verketteter {V}org\"{a}nge.
\newblock {\em Z. Angewandte Math. Mech.}, 3:279--289, 1923.

\bibitem{fri}
B.~Friedman.
\newblock A simple urn model.
\newblock {\em Communications on Pure and Applied Mathematics}, 2(1):59--70,
  1949.

\bibitem{hall-1980}
P.~Hall and C.~C. Heyde.
\newblock {\em Martingale limit theory and its application}.
\newblock Academic Press, Inc. [Harcourt Brace Jovanovich, Publishers], New
  York-London, 1980.
\newblock Probability and Mathematical Statistics.

\bibitem{Higueras2006}
I.~Higueras, J.~Moler, F.~Plo, and M.~San~Miguel.
\newblock Central limit theorems for generalized {P}{\'o}lya urn models.
\newblock {\em Journal of Applied Probability}, 43(4):938--951, 2006.

\bibitem{hu-ros}
F.~Hu and W.~F. Rosenberger.
\newblock {\em The Theory of Response-Adaptive Randomization in Clinical
  Trials}.
\newblock Wiley, Hoboken, 2006.

\bibitem{Idriss2018}
S.~Idriss and N.~Lasmar.
\newblock {Limit Theorems for Stochastic Approximations Algorithms With
  Application to General Urn Models}.
\newblock Hal-01726014, 2018.

\bibitem{Johnson2004}
N.~Johnson, S.~Kotz, and H.~Mahmoud.
\newblock {P{\'o}lya-Type Urn Models with Multiple Drawings}.
\newblock {\em J. Iran. Stat. Soc.}, 3(2):165--173, 2004.

\bibitem{kuba2016classification}
M.~Kuba.
\newblock Classification of urn models with multiple drawings.
\newblock Preprint Arxiv 1612.04354, 2016.

\bibitem{KubaMahmoud-balanced-affine-2017}
M.~Kuba and H.~M. Mahmoud.
\newblock Two-color balanced affine urn models with multiple drawings.
\newblock {\em Adv. in Appl. Math.}, 90:1--26, Sept. 2017.

\bibitem{mailler}
N.~Lasmar, C.~Mailler, and O.~Selmi.
\newblock {Multiple drawing multi-colour urns by stochastic approximation}.
\newblock {\em J. Appl. Probab.}, 55(1):254--281, 2018.

\bibitem{mahmoud_2013_multisets}
H.~M. Mahmoud.
\newblock Drawing multisets of balls from tenable balanced linear urns.
\newblock {\em Probability in the Engineering and Informational Sciences},
  27(2):147--162, 2013.

\bibitem{mayflournoy}
C.~May and N.~Flournoy.
\newblock Asymptotics in response-adaptive designs generated by a two-color,
  randomly reinforced urn.
\newblock {\em Ann. Statist.}, 37(2):1058--1078, 04 2009.

\bibitem{pemantle-volkov-1999}
R.~Pemantle and S.~Volkov.
\newblock Vertex-reinforced random walk on $\mathbf{Z}$ has finite range.
\newblock {\em Ann. Probab.}, 27(3):1368--1388, July 1999.

\bibitem{ros}
W.~F. Rosenberger.
\newblock Randomized urn models and sequential design.
\newblock {\em Sequential Anal.}, 21:1--28, 2002.

\bibitem{ros-lac}
W.~F. Rosenberger and J.~M. Lachin.
\newblock {\em Randomization in Clinical Trials}.
\newblock Wiley, New York, 2002.

\end{thebibliography}

\end{document}